\documentclass[11pt]{article}
\usepackage{epsfig}
\newtheorem{thm}{Theorem}

\font\msbm=msbm10 scaled 1200
\font\msbmscript=msbm8
\font\msbmscriptscript=msbm6
\newfam\Bbbfam
\textfont\Bbbfam=\msbm
\scriptfont\Bbbfam=\msbmscript
\scriptscriptfont\Bbbfam=\msbmscriptscript
\def\Bbb#1{{\fam\Bbbfam#1}}

\makeatletter
\@addtoreset{equation}{section}
\def\thethm{\arabic{thm}\protect\@blinkpoint}
\def\@blinkpoint{.}
\let\@blinkref=\ref
\def\ref#1{{\def\@blinkpoint{}\@blinkref{#1}}}
\let\@afterindentfalse\@afterindenttrue
\makeatother
\makeatletter
\long\def\@makecaption#1#2{%
   \vskip 10\p@
   \setbox\@tempboxa\hbox{#1. \ #2}%
   \ifdim \wd\@tempboxa >\hsize
       #1. \ #2\par
     \else
       \hbox to\hsize{\hfil\box\@tempboxa\hfil}%
   \fi}

\makeatother
\def\halmos{\hfill\rule{6pt}{6pt}}
\begin{document}
\thispagestyle{empty}
\medskip 
\begin{center}
  {\large\bf On the Set of Uniform  Convergence for the Last
  Intermediate Row of the  Pad\'e Table}
\footnote{This work was supported by
Russian Foundation for Basic Research under grant 04-01-96006.}

\medskip 

Victor M. Adukov \\
{ \it Department of Differential Equations and Dynamical Systems,\\
Southern Ural State University, Lenin avenue 76, 454080 Chelyabinsk, Russia\\
E-mail address: avm@susu.ac.ru} 
\end{center}
\medskip

{\small \it  
\noindent{\bf Abstract}

Let $a(z)$ be a meromorphic function having in the disk $|z|<R$ precisely
$\lambda$ poles. In this work  for the $(\lambda-1)$th row of the 
Pad\'e table of $a(z)$ the set of uniform convergence is explicitly 
obtained. The present note is a supplement to the previous work of the author ({\em J. Approx. Theory},
{\bf 123}(2003), 160-207).}
\medskip

In the theory of uniform convergence of the Pad\'e approximants the
principal question is if the presence of limit points of poles
for a sequence of the Pad\'e approximants in the disk ${\cal D}$ is the unique
obstraction for the unuform convergence of the sequence on compact subsets
of ${\cal D}$. For the diagonal sequence the affirmative answer has been given
(under some normality conditions) by A.A. Gon\v{c}ar~\cite{Gonchar}.

In the paper~\cite{Adukov} it was found all limit points of poles of the Pad\'e 
approximants for the row known as the last intermediate row of the 
Pad\'e table for  a meromorphic function. In the present note we are going
to show that limit points of poles for this row are also the unique
obstraction for the uniform convergence. Thus taking into account the results
of \cite{Adukov} we can obtain the set of uniform convergence for the last 
intermediate row.

Recall relevant definitions and statements. Let $a(z)$  be a function which is
meromorphic in the disk  ${\cal D}_R=\left\{z\in\Bbb C\Bigl||z|< R\right\}$
and analytic at the origin. Let $z_1,\ldots,z_{\ell}$ 
be its distinct poles of multiplicities 
$s_1,\ldots,s_{\ell}$, respectively, and let
$\lambda=s_1+\ldots+s_{\ell}$ be the number of its 
poles in the disk ${\cal D}_R$. Suppose that 
$
\rho\equiv |z_1| =\ldots=|z_{\mu}|>|z_{\mu +1}|\ge 
\ldots \ge |z_{\ell}|.
$

If 
$m=\sum_{j=1}^{\ell}s_j$ or $m=\sum_{j=\mu+1}^{\ell}s_j$,  
then, by de Montessus's theorem, the Pad\'e approximants 
$\pi_{n,m}(z)$ converge to
$a(z)$, as $n\to\infty$, uniformly on compact subsets of the set
${\cal D}_R\setminus\left\{z_1,\ldots,z_{\ell}\right\}$ or
${\cal D}_\rho\setminus\left\{z_{\mu+1},\ldots,z_{\ell}\right\}$, 
respectively.
The row of the Pad\'e table with the number $m$ satisfying
the inequalities 
$$
\sum_{j=\mu+1}^{\ell}s_j<m<\sum_{j=1}^{\ell}s_j
$$
is said to be an {\em intermediate row}. Sufficient conditions
for the convergence of the whole intermediate row has been got in
\cite{Sidi}. 

For the row with the number 
$m=\lambda-1=\sum_{j=1}^{\ell}s_j-1$ (the last intermediate row)
 the asymptotic behavior of denominators 
$Q_{n,\lambda-1}(z)$ for the Pad\'e approximants 
$\pi_{n,\lambda-1}(z)$ and all limit points of poles of 
$\pi_{n,\lambda-1}(z)$ are known. Let us describe these 
results~\cite{Adukov}.
The poles $z_1,\ldots,z_{\mu}$ of the maximal modulus
we will order in such a way that
$s_1\ge\ldots\ge s_{\mu}$.
Among the poles $z_1,\ldots,z_{\mu}$ 
we select the poles 
$z_1,\ldots,z_{\nu}$ ($1\le \nu \le \mu\le \ell$)
that have the maximal multiplicity:
$
s_1=\ldots=
s_{\nu}>s_{\nu+1}\ge\ldots\ge s_{\mu}.
$
The poles  $z_1,\ldots,z_{\nu}$ will be called 
{\em dominant poles} of $a(z)$.
Let
$
z_1 = \rho e^{2\pi i\Theta_1}, \ldots,
z_{\nu} = \rho e^{2\pi i\Theta_{\nu}}.
$
Consider the point $\xi =\left (e^{2\pi i \Theta_1}, \ldots,
e^{2\pi i \Theta_{\nu}}\right)$ belonging to the torus $\Bbb T^{\nu}$.
The torus is a compact Abelian group. Denote by $\Bbb F$ the closure 
in $\Bbb T^{\nu}$ of the cyclic group $\left\{\xi^n\right\}_{n\in \Bbb Z}$
with the generator $\xi$. $\Bbb F$ is a monothetic subgroup of the 
torus $\Bbb T^{\nu}$. Let $r+1$ be the rank over the field of rational
numbers $\Bbb Q$ of the system of the real numbers 
$\Theta_0=1,\Theta_1, \ldots,\Theta_{\nu}$. If $r=\nu$, 
then $\Bbb F=\Bbb T^{\nu}$. For $0\le r<\nu$ the group $\Bbb F$ is
isomorphic to ${\Bbb Z}_{\sigma}\times {\Bbb T}^r$. In this case 
the group $\Bbb F$ can be explicitly found if the matrix of linear
relations between the arguments $\Theta_1, \ldots,\Theta_{\nu}$
is known (see~\cite{Adukov}, Theorem~2.1). 
In the problem under consideration the group $\Bbb F$ plays a
significant role. The limits of all convergent subsequences of
the sequence $Q_{n,\lambda-1}(z)$ form a family of polynomials which is
parametrized by $\Bbb F$.

Let  $A_j$ be the coefficient of $(z-z_j)^{-s_j}$ in the Laurent series  
in a neighborhood of the pole  $z=z_j$ for the function $a(z)$.
Put
$$
C_j=\frac{1}{(s_j-1)!z_j^{s_j-1}D^2_j(z_j)A_j},
\leqno (1)
$$
where
$
D_j(z) = \frac{D(z)}{(z-z_j)^{s_j}}, \ \ D(z)=(z-z_1)^{s_1}\cdots (z-z_\ell)^{s_\ell}.
$

Then the set of the limit points of poles of the sequence
$\pi_{n,\lambda-1}(z)$, as $n\to\infty$, consists of the poles 
$z_1,\ldots,z_\ell$ (the multiplicity of the dominant poles $z_1,\ldots,z_{\nu}$ 
is less by 1), and the set ${\cal N}_{\Bbb F}$ of the zeros
of polynomials from the family 
$$
\omega(z,\tau)=\sum_{j=1}^{\nu}C_j\Delta_j(z)\tau_j,\ \ 
\tau=\left(\tau_1,\ldots,\tau_\nu\right)\in\Bbb F.
$$
Here $\Delta_j(z)= \frac{\Delta(z)}{z-z_j}, \ \ \Delta(z)=(z-z_1)\cdots (z-z_{\nu}).$

The set ${\cal N}_{\Bbb F}$ is a closed set and $z_1,\ldots,z_\nu$ does not belong
to it. Let ${\cal N}$ be the set of complex points $z$ satisfying the 
following system of inequalities
$$
2\bigl|C_j\Delta_j(z)\bigr| \le
\sum_{k=1}^{\nu} 
\bigl|C_k\Delta_k(z)\bigr|,
\ \ j=1,\ldots,\nu.
$$
Then ${\cal N}_{\Bbb F}\subseteq{\cal N}$. Moreover,
if $r=\nu$, i.e. the numbers $\Theta_0=1,\Theta_1, \ldots,\Theta_{\nu}$ are linearly 
independent over $\Bbb Q$, then ${\cal N}_{\Bbb F}={\cal N}$ (see~\cite{Adukov},
Theorem 2.6).

Let us delete the set ${\cal N}_{\Bbb F}\cap {\cal D}_\rho$ 
(${\cal N}\cap{\cal D}_\rho$) and  the poles 
$z_{\mu+1},\ldots,z_\ell$ from ${\cal D}_\rho$. The open set obtained we denote by 
${\Bbb U}_{\Bbb F}$ (${\Bbb U}$).
Evidently, ${\Bbb U}\subseteq {\Bbb U}_{\Bbb F}$. Moreover, for every
dominant pole $z_j$, $j=1,\ldots,\nu$, there exists a neighborhood $U_j$
such that $U_j\cap {\cal D}_\rho \subseteq \Bbb U$. Our aim is to prove that
the sequence $\pi_{n,\lambda-1}(z)$ uniformly converges on compact
subsets of ${\Bbb U}_{\Bbb F}$ to $a(z)$ as $n\to\infty$.

To do this we will need the simplified version of the Vitali theorem.
Recall that a sequence of analytic functions $f_n(z)$ is called 
{\em compact} in an open set $G$ if from each subsequence $f_{n_k}(z)$
we can select a subsequence $f_{n_{k_i}}(z)$ that uniformly converges 
on compact subsets of $G$. 
\begin{thm}[Vitali]
Let a sequence $f_n(z)$ be compact in an open set $G$. If all uniformly 
convergent on compact subsets of $G$ subsequences $f_{n_{k_i}}(z)$ have 
the same limit function $f(z)$, then the sequence $f_{n}(z)$
uniformly converges on compact subsets of $G$ to $f(z)$.
\halmos
\end{thm}

In contrast to the standart version of Vitali's theorem in Theorem 1 we
require that all subsequences $f_{n_{k_i}}(z)$ converge to 
the same limit function $f(z)$. The proof of the Theorem 1 coincides with 
the proof of the second part of Vitali's theorem (see \cite{Mark}, p. 371).

\begin{thm}
The sequence $\pi_{n,\lambda-1}(z)$ uniformly converges on compact subsets of 
$\Bbb U_{\Bbb F}$ to $a(z)$.
\end{thm}
{\bf Proof.} Let us prove that the sequence $\pi_{n,\lambda-1}(z)$ is compact
in the set $\Bbb U_{\Bbb F}$. Take any sequence of natural numbers
$n_1, n_2,\ldots,n_k,\ldots $\ .  From the set of points
$\xi^{n_k+\lambda} =\left (e^{2\pi i(n_k+\lambda) \Theta_1}, \ldots,
e^{2\pi i (n_k+\lambda)\Theta_{\nu}}\right)\in\Bbb T^{\nu}$ we can
select a subsequence $\xi^{n_{k_i}+\lambda}$ that converges to some 
point $\tau_0\in\Bbb T^{\nu}$.
It follows from the definition of the group ${\Bbb F}$ that 
$\tau_0\in\Bbb F$. Let us denote by $\Lambda_{\tau_0}$ the 
subsequence $\left\{n_{k_i}+\lambda\right\}_{i=1}^\infty$.
Let $K_{\tau_0}$ be any compact subset of the dick ${\cal D}_\rho$ such
that all zeros of the polynomial $\omega(z,\tau_0)$ and the poles
$z_{\mu+1},\ldots,z_\ell$ lie outside $K_{\tau_0}$. By Theorem~2.3
in~\cite{Adukov} the subseqence 
$\left\{\pi_{n,\lambda-1}\right\}_{n\in\Lambda_{\tau_0}-\lambda}$
uniformly converges to $a(z)$ on $K_{\tau_0}$. This means that 
the sequence $\pi_{n_{k_i},\lambda-1}(z)$ uniformly converges to
$a(z)$ on compact subsets of $\Bbb U_{\Bbb F}$, i.e. the 
sequence $\pi_{n,\lambda-1}(z)$ is compact in the 
set $\Bbb U_{\Bbb F}$. To conclude the proof of the theorem, it
remains to apply Theorem 1. \halmos

Here we give an example of the determination of the set 
$\Bbb U_{\Bbb F}$. For calculations we used Maple6.

\vspace{2cm}

\begin{figure}[h]  \label{Fig1} 
\epsfig{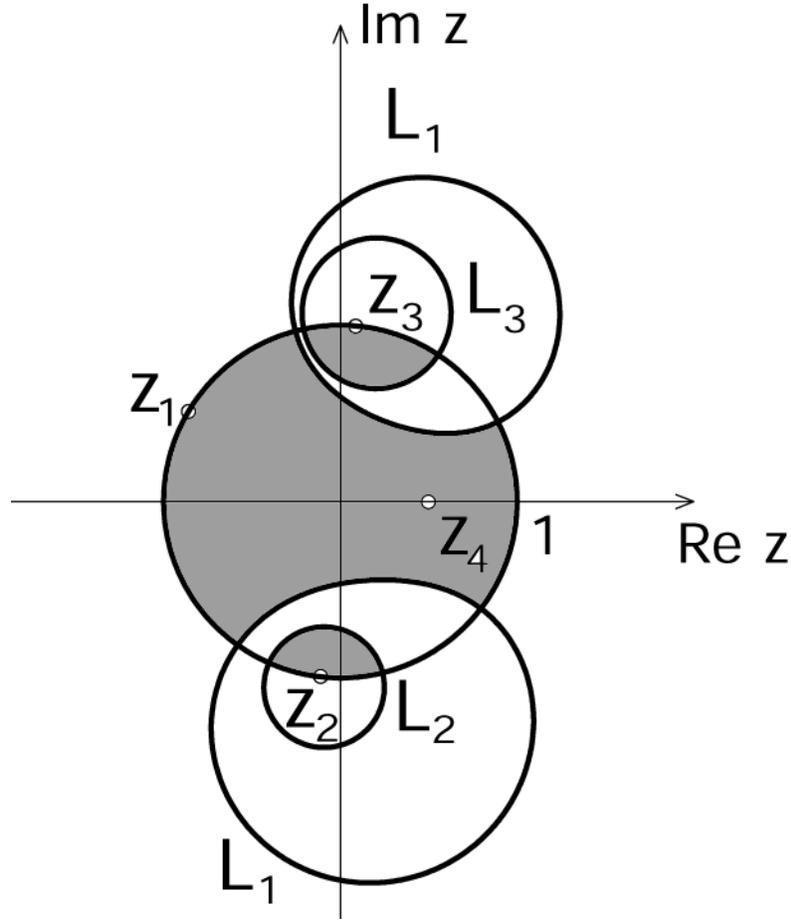}
\caption{The example of the set of uniform convergence}
\end{figure}

{\noindent\bf Example.} Let $b(z)$ be any function which is analytic
in the disk ${\cal D}_R$, $R>1$. Let
$$
r(z)=
\frac{z^2+z}{(z-e^{2\pi\imath\sqrt{2}})
(z-e^{2\pi\imath\sqrt{3}})
(z-e^{2\pi\imath\sqrt{5}})
(z-1/2)},
$$
and $a(z)=b(z)+r(z)$.

The dominant poles of the meromorphic function $a(z)$ are
$z_1=e^{2\pi\imath\sqrt{2}}$, $z_2=e^{2\pi\imath\sqrt{3}}$,
$z_3=e^{2\pi\imath\sqrt{5}}$. 
Since $\Theta_0=1$, $\Theta_1=\sqrt{2}$, $\Theta_2=\sqrt{3}$,
$\Theta_3=\sqrt{5}$ are linearly independent over $\Bbb Q$,
we have ${\cal N}_{\Bbb F}={\cal N}$ and  $\Bbb U_{\Bbb F}=\Bbb U$.
By formula (1) we obtain
$
{C_{1}} = 0.70400 + 0.17095\,i, \ \ 
{C_{2}} = 0.07853 + 0.17437\,i, \ \ 
{C_{3}} = 0.29275 + 0.04487\,i.
$

The boundary of the set ${\cal N}$ consists of the lines $L_1, L_2,
L_3$, where $L_j$ is given by the following equation:
$$
2\bigl|C_j\Delta_j(z)\bigr| =
\sum_{k=1}^{3} 
\bigl|C_k\Delta_k(z)\bigr|,
\ \ j=1, 2, 3.
$$
These lines are shown in Figure 1. $L_1$ is a reducible curve
consisting of  two connected components. The set 
$\Bbb U_{\Bbb F}$ is represented by the dark-shaded area.

If we have known nothing about an arithmetic nature of the 
dominant poles, we can claim that the sequence 
$\pi_{n,\lambda-1}(z)$ uniformly converges on compact subsets of
$\Bbb U$. At last, if the Laurent coefficients $A_j$ of $a(z)$ are
also unknown, we can only assert that for each dominant pole $z_j$
there exists a neighborhood $U_j$ such that $\pi_{n,\lambda-1}(z)$    
uniformly converges on $U_j\cup {\cal D}_{\rho}.$

It remains to answer the following principle question.
Are the limit points of poles of the Pade approximants really the
obstruction for the uniform convergence? This will be the 
subject of another paper.


\end{document}